\documentstyle{amsppt}
\magnification 1200 \TagsOnRight
\parskip=\medskipamount
\def\ve{\varepsilon}
\def\vp{\varphi}
\def\arrowk{^\to{\kern -6pt\topsmash k}}
\def\arrowK{^{^\to}{\kern -9pt\topsmash K}}
\def\arrowt{^\to{\kern -6pt\topsmash t}}
\def\arrowr{^\to{\kern-6pt\topsmash r}}
\def\arrowvp{^\to{\kern -8pt\topsmash\vp}}
\def\tk{\tilde{\kern 1 pt\topsmash k}}

\hsize = 6.2true in \vsize=8.2 true in \NoRunningHeads \topmatter
\title
{\smc\bf On the size of k-fold Sum and Product Sets of Integers}
\endtitle
\author
{Jean Bourgain }\quad\quad {Mei-Chu Chang}
\endauthor
\address{Institute for Advanced Study, Olden Lane,
Princeton, N.J. 08540, U.S.A.}
\endaddress
\email{bourgain\@math.ias.edu}
\endemail
\address{Mathematics Department, University of California,
Riverside, CA 92521}
\endaddress
\email {mcc\@math.ucr.edu}
\endemail
\endtopmatter
\document

\noindent {\bf Section 1. Statement of the result and outline of
argument}

We prove the following
\medskip
\proclaim {Theorem} For all $b\in\Bbb Z_+$, there is
$k=k(b)\in\Bbb Z_+$ such that if $A\subset\Bbb Z$ is any finite
set, with $|A|=N\geq 2$, then either
$$
|kA|= |\underbrace{A+\cdots +A}_k| > N^b,\tag 1.1
$$
or
$$
|A^{(k)}|=|\underbrace{A\times\cdots\times A}_k|>N^b.\tag 1.2
$$
\endproclaim

This result is one more contribution to a line of research
initiated in the seminal paper [E-S] of Erd\"os and Szemer\'edi on
sum and product sets. They conjectured that if $A \subset \Bbb Z$,
with $|A| =N$, then
$$
|2A|+|A^{(2)}| > c_{\varepsilon} N^{2-\varepsilon}, \quad \text
{for all } \varepsilon > 0, \tag i
$$
and more generally, for $k\geq 2$ an integer,
$$
|kA|+|A^{(k)}| > c_{\varepsilon} N^{k-\varepsilon}, \quad \text
{for all } \varepsilon > 0, \tag ii
$$

Already (i) is open. Recent advances were achieved by G. Elekes
[E] and J. Solymosi [So] using the Szemer\'edi-Trotter theorem in
incidence geometry. It was shown in [E] that
$$
|2A| +|A^{(2)}|> cN^{5 \over 4},
$$
and in [So] a small further improvement
$$
|2A| + |A^{(2)}|> cN^{{14 \over 11}-\varepsilon}.
$$

Again, using the Szemer\'edi-Trotter theorem,
Elekes-Nathanson-Ruzsa showed that
$$
|kA|\cdot|A^{(k)}|>c|A|^{3-2^{1-k}}.
$$
(See [E-N-R].)  The theorem proved in this paper answers
affirmatively a problem posed in their paper and provides further
progress towards (ii). This same issue was also brought up
independently by S. Konjagin (private communication) and was
motivated by questions related to Gauss-sums (see [B-K]). Some
further refinements of the statement of the theorem appear in
Section 6 (Remarks) at the end of the paper.

The general strategy of our proof bares resemblance with [Ch] in
several ways. Thus we assume $|A^{(k)}|$ `small' and prove that
then $|kA|$ has to be large. However, `smallness' of $|A\cdot A|$
in [Ch] is the assumption
$$
|A\cdot A|<K|A|\tag 1.4
$$
with $K$ a constant (a condition much too restrictive for our
purpose).

If (1.4) holds, it is shown in [Ch] that
$$|A+A|>c(K)|A|^2\tag 1.5
$$
and more generally
$$
|hA|> c(K, h)|A|^h.\tag 1.6
$$
Let us briefly recall the approach.

Consider the map given by prime factorization
$$
\Bbb Z_+\longrightarrow \Cal R= \prod_p\Bbb Z_{\geq 0}
$$
$$
n=\prod_p p^{\alpha_p}\longrightarrow \alpha =(\alpha _p)_p
$$
where $p$ runs in the set $\Cal P$ of primes.

The set $A$ is mapped to $\Cal A\subset\Cal R$ satisfying by (1.4)
$$
|2\Cal A|<K|\Cal A|.\tag 1.4'
$$
Freiman's lemma implies then that $\dim \Cal A< K$ (where `dim'
refers to the dimension of the smallest vector space containing
$\Cal A$). Hence there is a subset $I\subset\Cal P, |I|<K$ such
that the restriction $\pi_I$ is one-to-one restricted to $\Cal A$.
Harmonic analysis implies then that
$$
\lambda_q(A)<(Cq)^K\tag 1.7
$$
for an absolute constant $C$, and for all $q>2$. By
$\lambda_q(A)$, we mean the $\Lambda_q$-constant of the finite set
$A\subset\Bbb Z$, defined by
$$
\lambda_q(A)=\max \bigg\Vert\sum_{n\in A} c_n \ e^{2\pi i n\theta}\bigg\Vert_{L^q(\Bbb T)}\tag 1.8
$$
where $\Bbb T=\Bbb R/\Bbb Z$ and the max is taken over all
sequences $(c_n)_{n\in A}$ with $\big(\sum c_n^2\big)^{1/2}\leq
1$. See [Ru] and [Ki], for more details. As is [Ch], we will also
make here crucial use of certain moment inequalities involving
$\lambda _q$-constant of certain specific sets of integers. More
specifically, we will use the following general inequality, from
which (1.7) is derived.

\proclaim
{Proposition 1} (see [Ch]).

Let $p_1, \ldots, p_k$ be distinct primes and associate to each
$\alpha =(\alpha_1, \ldots, \alpha_k) \in(\Bbb Z_{\geq 0})^k$ a
trigonometric polynomial $F_\alpha$ on $\Bbb T$ such that
$$
(n, p)=1, \text { for all  } n\in\text{\,supp\,}\widehat F_\alpha,
\text {  and for all } p \in \Cal P_0.
$$
Then, for any moment $q\geq 2$
$$
\bigg\Vert \sum_\alpha F_\alpha (p_1^{\alpha_1} \cdots p_k^{\alpha_k}\theta)\bigg\Vert_q < (Cq)^k \bigg(\sum\Vert
F_\alpha\Vert^2_q\bigg)^{1/2}.\tag 1.9
$$
\endproclaim

Thus (1.7) follows from (1.9) taking $F_\alpha(\theta) =e^{2\pi
i\theta}$ and $\{p_1, \cdots, p_k\}=I\subset\Cal P$.

Denoting for $h\geq 2$
$$
r_h(n; A)=|\{(x_1,\ldots, x_h)\in A^h|n=x_1+\cdots+x_h\}|
$$
A simple application of Parseval's identity gives
$$
\sum_{n\in hA} r_h(n; A)^2\leq \lambda_{2h}(A)^{2h}\cdot |A|^h
$$
and using Cauchy-Schwartz inequality on $\sum_{n\in hA} r_h(n;
A)$, it follows that
$$
|hA|\geq\frac{|A|^h}{\lambda_{2h}(A)^{2h}}.\tag 1.10
$$
Thus we obtain (1.6) with
$$
c(K, h)> (Ch)^{-2hK}.
$$
Obviously, this statement has no interest unless $K\ll \log|A|$.

The main point in what follows is to be able to carry some of the
preceding analysis under a much weaker assumption $K<|A|^\ve$,
$\ve$ small. We will prove the following statement

\proclaim {Proposition 2} Given $\gamma>0$ and $q>2$, there is a
constant $\Lambda=\Lambda(\gamma, q)$ such that if $A\subset\Bbb
Z$ is a finite set, $|A|=N, |A\cdot A|<KN$, then
$$
\lambda_q(A)<K^\Lambda N^\gamma\tag 1.11 $$ \endproclaim Thus
fixing $q$, Proposition 2 provides already nontrivial information
assuming $K <N^{\delta}$, with $\delta >0$ sufficiently small.

Assuming Proposition 2, let us derive the Theorem. We may assume
that $A \subset \Bbb Z_+$ to simplify the situation.

Fix $b$ and assume (1.2) fails for some large $k=2^\ell$ (to be
specified). Hence, passing to $\Cal A$
$$
|k\Cal A| <N^b
$$
$$
\frac{|2^\ell \Cal A|}{|2^{\ell-1}\Cal A|} \ \frac{|2^{\ell-1}\Cal
A|}{|2^{\ell-2}\Cal A|} \cdots \frac{|2\Cal A|}{|\Cal A|} <
N^{b-1}
$$
and we may find $k_0=2^{\ell_0}$ such that
$$
\frac{|2k_0\Cal A|}{|k_0\Cal A|}< N^{\frac{b-1}\ell}.\tag 1.12
$$
Denote $\Cal B=k_0\Cal A\subset\Cal R$ and $B=A^{(k_0)}$,
 the
corresponding subset of $\Bbb Z_+$. Thus by (1.12)
$$
|B\cdot B|<N^{\frac{b-1}\ell}|B|.\tag 1.13
$$
Apply Proposition 2 to the set $B$, $|B|\equiv N_0$,
$K=N^{\frac{b-1}\ell}$ with $\tau, \gamma$ specified later.

Hence from (1.11)
$$
\lambda_q(A)\leq\lambda_q(B)<N^{\frac{b-1}\ell \Lambda}
N_0^{\gamma}<N^{\frac{b-1}\ell \Lambda +b\gamma} \tag 1.14
$$

Taking $q=2h$, (1.10) and (1.14) imply
$$
|hA|> N^{(1-2\frac{b-1}{\ell}\Lambda -2b\gamma)h} \tag 1.15
$$
Take $h=2b<k, \gamma =\frac {1}{100b}.$ Recall that $\Lambda
=\Lambda (\gamma, q)$, hence $\Lambda =\Lambda (b).$ Take $\ell =
100 \,b \,\Lambda (b)$, so that $k=2^{\ell}\equiv k(b)$.
Inequality (1.15) then clearly implies that
$$
|kA|> N^b
$$
This proves the Theorem. \medskip Returning to Proposition 2, it
will suffice to prove the following weaker version

\proclaim {Proposition 2'} Given $\gamma>0, \tau >0$ and $q>2$,
and $A$ as in Proposition 2, there is a subset $A' \subset A$
satisfying
$$
|A'|>N^{1-\tau}\tag 1.16 $$ $$ \lambda _q(A')<
K^{\Lambda}N^{\gamma},\tag 1.17
$$ where $\Lambda =\Lambda (\tau,\gamma,q). $
\endproclaim
\medskip
\noindent {\bf Proof of Proposition 2 assuming Proposition 2'.}

Denoting $\chi$ the indicator function, one has obviously $$
\sum_{z\in \frac {A}{A'}} \chi_{z A'} \geq |A'|\, \chi _A \tag
1.18
$$

Let $A'$ be the subset obtained in Proposition 2'. Then (1.18) is
easily seen to imply
$$
\align |A'| \,\lambda _q(A)\,& \leq \sum_{z\in\frac{A}{A'}}
\lambda_q(z
 A')\\
&= |\frac {A}{A'}|\,
\,\lambda_q(A')\\
&\leq |\frac {A}{A}\, |\,K^{\Lambda}N^{\gamma}.\tag 1.19
\endalign
$$

If $\Cal A \subset \Cal R$ is the set introduced before,
application of Ruzsa's inequality on sum-difference sets gives
$$
|\frac {A}{A}|=| \Cal A - \Cal A| \leq K^2\, |\Cal A| = K^2N.\tag
1.20
$$

Thus, by (1.16), (1.19) and (1.20), we have
$$
\lambda _q (A) \leq K^{\Lambda +2}N^{\tau + \gamma}, \tag 1.21
$$
where $\Lambda =\Lambda (\tau, \gamma, q)$. Replacing $\gamma$ by
$\frac {\gamma}{2}$ and $\tau = \frac {\gamma}{2}$, (1.11)
follows.
\medskip
The remainder of the paper is the proof of Proposition 2' which
will be rather tedious although elementary. The key statement is
Proposition 3 below. This more technical result involves "graphs".
We were unable to carry out our analysis by passing simply to
subsets. As will be clear later on, the use of graphs allows
indeed more flexibility in various constructions involving certain
"regularizations".  Trying to achieve them using subsets of $A$ is
less economical and did not seem to provide us with the desired
conclusions. These considerations are particularly relevant to
Lemma 3.2 below (which is the base of the multi-scale analysis)and
the difficulties encountered with its proof.

\bigskip
\noindent {\bf Section 2. Reduction to Proposition 3 and
preliminary estimates}

Given subsets $A_1, A_2$ in an Abelian group and a symmetric graph $\Cal G\subset A_1\times A_2$, define
$$
A_1\underset {\Cal G}\to + A_2 =\{x_1+x_2|(x_1, x_2)\in\Cal G\}.
$$
For technical reason, we will prove a more general (and
complicated) version of Proposition 2 (we consider subsets of
$\Cal R$, hence work in the additive setting).

\proclaim {Proposition 3} Given $\tau, \gamma>0$ and $q$, there is
a constant $\Lambda=\Lambda(\tau, \gamma, q)$ such that the
following holds.

\noindent Let $\Cal P_0\subset\Cal P$ be a set of primes, and let
$\Cal R_0=\prod_{p\in\Cal P_0}\Bbb Z_{\geq 0}$. Let $A_1$, $A_2
\subset \Cal R_0$ be finite with $|A_i| = N_i$ and let
$$
N=N_1N_2. $$
 If $\Cal G\subset A_1\times A_2$ with
$$
|\Cal G|>\delta N,
$$
then there is $\Cal G'\subset\Cal G$ satisfying
$$
|\Cal G'|>\delta^{\Lambda\log\log N}N^{1-\tau}\tag 2.1
$$
and, moreover, for all $x\in\Cal R_0$, the set $\Cal G'(x)=\{x'
\mid (x,x')\in \Cal G', \text { or } (x', x)\in \Cal G' \}$ has
the following property:

\noindent If $(F_\alpha)_{\alpha\in \Cal G'(x)}$ are arbitrary
trigonometric polynomials such that $(n, p)=1$ for all
$n\in\text{\,supp\,}\widehat F_\alpha$ and for all $p \in \Cal
P_0$, then
$$ \bigg\Vert \sum_\alpha F_\alpha \bigg(\prod_{\Cal P_0}
p^{\alpha_p}\theta\bigg)\bigg\Vert_q \leq K(\Cal G)^{\Lambda}
N^\gamma \bigg(\sum_\alpha \Vert F_\alpha\Vert^2_q\bigg)^{1/2}\tag
2.2
$$
holds, where we denote
$$
K(\Cal G)=\frac{|A_1\underset{\Cal G}\to+ A_2|}{\sqrt {N_1N_2}}.
$$
\endproclaim

\medskip
\noindent {\it Remark.} If $A_1 =A_2=A$ and $\Cal G = A \times A$,
then $K(\Cal G)$ is the doubling constant of $A$.
\medskip
We will use the following

\noindent {\bf Notation.} $M \sim N$ means that there are
constants $c$ and $d$ such that $dN<M<cN$.
\medskip

\noindent {\bf Proposition 3 implies Proposition 2'.} Take $A_1=
A_2=\Cal A\subset\Cal R$, $\Cal P_0=\Cal P$, and the full graph
$\Cal G=A_1\times A_2$. (Note that $K(\Cal G) =K$, the doubling
constant in Proposition 2.) Since $\Cal G' \subset \Cal A \times
\Cal A $ satisfies, by (2.1) $$ |\Cal G'|> N^{2-3\tau}, $$ there
is $x\in \Cal A$ such that $|\Cal G'(x)|>N^{1-3\tau}$. Let
$A'\subset A$ be the set corresponding to $\Cal G'(x)$. Thus
$|A'|>N^{1-3\tau}$ and, taking
$F_{\alpha}(\theta)=c_{\alpha}e^{2\pi i \theta}$, $c_{\alpha} \in
\Bbb R$, statement (2.2) clearly implies (1.17)
\medskip

The proof of Proposition 3 will proceed in several stages. In this
process, we consider pairs of functions
$$
\phi(N, \delta, K), \psi(N, \delta, K)
$$
such that under the assumptions of Proposition 3, there is $\Cal G'\subset\Cal G$ satisfying
$$
|\Cal G'|>\phi (N, \delta, K(\Cal G))\qquad (N=N_1N_2)\tag 2.3
$$
and (2.2) holds in the form
$$
\bigg\Vert\sum_\alpha F_\alpha \bigg(\prod p^{\alpha_p}\theta\bigg)\bigg\Vert_q \leq \psi(N, \delta, K(\Cal G))
\bigg(\sum_\alpha \Vert F_\alpha\Vert^2_q\bigg)^{1/2}.\tag 2.4
$$
(We assume $q$ fixed.) We call such a pair of functions {\it
admissible}.

The strategy then consists in getting better and better bounds on
the functions $\phi(N, \delta, K), \psi(N, \delta, K)$ and
eventually prove Proposition 3. Let us specify a first pair of
admissible functions $\phi, \psi$. (see (2.13), (2.14))

Assume $N_1\geq N_2$.

Obviously $|A_1\underset{\Cal G}\to+ A_2|\geq \frac {|\Cal
G|}{N_2}\geq \delta{N_1}$, hence
$$
\delta N_1 < K(\Cal G)(N_1N_2)^{1/2}. $$ Namely,
$$
N_2>\bigg(\frac{\delta} {K(\Cal G)}\bigg)^2 N_1.\tag 2.5
$$

We may assume that $\Cal G\subset A\times A$ is symmetric, where
$A=A_1\cup A_2,$ with $ |A|\sim N_1$. Let $$K=K(\Cal G),$$ we have
thus
$$
\align
&|\Cal G|> \delta N_1 N_2>\frac{\delta^3}{K^2}N_1^2 =\delta_1N_1^2\\
&|A\underset{\Cal G}\to +  A|=|A_1\underset{\Cal G} \to +A_2|<
KN_1,
\endalign
$$
where $$\delta _1 = \frac { \delta ^3}{K^2}.$$

 Applying Gower's
version of the Balog-Szemeredi theorem (with powerlike estimate),
see [Go], we may find  a subset $A'\subset A$ satisfying
$$
\align
&|A'|>\delta'N_1\tag 2.6\\
&|A'-A'|< K'N_1\tag 2.7\\
&|\Cal G\cap (A'\times A')|>\delta'N_1^2\tag 2.8
\endalign
$$
where
$$
\delta' >\bigg(\frac{\delta_1}K\bigg)^C>\bigg(\frac \delta
K\bigg)^C,
$$
$$
 K' <\bigg(\frac K\delta\bigg)^C\tag 2.9
$$
(here and in the sequel, notation $C$ as well as $c$ may refer to
different constants).

By (2.6), (2.7) and Freiman's lemma, the dimension of the vector
space spanned by $A'$ is less than $\frac{K'}{\delta'}$. Thus
there is a subset $I\subset \Cal P_0$ such that
$$
|I|<\frac{K'}{\delta'} \tag 2.10
$$
and the coordinate restriction $\pi_I:\Cal R_0\rightarrow
\prod_{p\in I}\Bbb Z_{\geq 0}$ is one-to-one when restricted to
$A'$.

Let $\Cal G'=\Cal G\cap(A'\times A')$ satisfying by (2.8), (2.9)
$$
|\Cal G'|>\bigg(\frac\delta K\bigg)^C N_1N_2.\tag 2.11
$$
Fix $x\in \Cal R_0$ and consider $\Cal G'(x)\subset A'$ and
trigonometric polynomials $F_\alpha, \alpha\in \Cal G'(x)$ as in
(2.2). Thus $(n, p)=1$, for all $n\in \text{\,supp\,}\widehat
F_\alpha,$ and for all $ p\in\Cal P_0$. It follows from the
preceding that $\alpha\in\Cal G'(x)$ is uniquely determined by
$\alpha'=\pi_I(\alpha)$. Therefore clearly
$$
F_\alpha\bigg(\prod_{p\in\Cal P_0}p^{\alpha_p}\theta\bigg)=F_{\alpha'}'\bigg(\prod_{p\in I} p^{\alpha_p}\theta\bigg)
$$
where $(n, p)=1$ for $n\in \text{\,supp\, } \widehat
F_{\alpha'}',$ and $ p\in I$.

Thus Proposition 1 and (2.10) imply
$$
\align
\bigg\Vert \sum_{\alpha\in\Cal G'(x) } F_\alpha (\prod p^{\alpha_p}\theta)\bigg\Vert_q
&=\bigg\Vert\sum_\alpha F_{\alpha'}' \bigg(\prod_{p\in I} p^{\alpha_p}\theta\bigg)\bigg\Vert_q\\
&\leq (Cq) ^{|I|} \bigg(\sum_\alpha\Vert F_{\alpha'}' \Vert_q^2\bigg)^{1/2}\\
&\leq (Cq)^{\frac{K'}{\delta'}} \bigg(\sum_\alpha \Vert F_\alpha\Vert^2_q\bigg)^{1/2}.\tag 2.12
\endalign
$$
Hence, (2,11), (2.12), and (2.9) provide the following pair of
admissible functions
$$
\align
\phi(N, \delta, K)&=\bigg(\frac\delta K\bigg)^CN\tag 2.13\\
\psi (N, \delta, K)&= \exp \bigg(\log q\cdot \bigg(\frac K\delta\bigg)^C\bigg)\tag 2.14
\endalign
$$
for some constant $C$.

Again the dependence of $\psi$ on $K$ is very poor, since it is a useless bound unless $K\ll \log N$.

The aim of what follows is to improve this dependence of $\psi$ on $K$.

\bigskip

\noindent {\bf Section 3. Proof of Proposition 3, Part I: The
Factorization}

The next statement is a recipe to convert pairs of admissible functions $\phi, \psi$.
We will always assume

$\phi, \psi$ are increasing in $N$

$\phi$ is increasing in $\delta$, decreasing in $K$

$\psi$ increases in $K$

\noindent
and
$$
\phi(N, \delta, K)\leq \frac NM \phi(M, \delta, K) \text { for } M\leq N.\tag 3.1
$$
\medskip
\proclaim
{Lemma 3.2} Let $\phi, \psi$ be admissible.
Define
$$
\align
&\tilde\phi (N, \delta, K)= \min \phi(N',\delta', K')\cdot \phi(N'', \delta'', K'')\tag 3.3
\\
&\tilde\psi(N, \delta, K)= Cq \max \psi (N', \delta', K')\cdot \psi(N'', \delta'', K'')\tag 3.4
\endalign
$$
where in (3.3), (3.4) the range of $N', N'', \delta', \delta'', K', K''$ are as follows
$$
\align
&N\geq N'N'' >N\bigg( \frac\delta{\log K}\bigg)^{40}\tag 3.5\\
&N'+N''< \bigg(\frac K\delta\bigg)^{20} N^{1/2}\tag 3.6\\
&\delta'\cdot\delta''>\bigg(\log\frac K\delta\bigg)^{-6}\delta\tag 3.7
\endalign
$$
$$
K'\cdot K'' < \delta^{-6}(\log K)^{20} K.\tag 3.8
$$
Then $\tilde \phi, \tilde\psi$ are also admissible.
\endproclaim

This lemma is an essential ingredient in the proof of Proposition
3. In its proof, the role of graphs will become apparent.
\medskip
Under the assumption of Proposition 3, we have $\Cal G\subset
A_1\times A_2,$ for $ A_i \subset \Cal R_0 =\prod_{p\in\Cal P_0}
\Bbb Z_{\geq 0},$ with $|A_i|=N_i$, and
$$
|\Cal G|> \delta N,
$$
$$|A_1\underset {\Cal G}\to +A_2|< K(\Cal G)\sqrt
N,
$$
where $N=N_1N_2$.
\medskip
The proof of Lemma 3.2 is in seven steps.
\medskip
\proclaim {Step 1} For $i= 1, 2$, we reduce $A_i$ to $A_i',$ with
$ |A_i'| =N_i'$ such that for any $B_i \subset A_i',$
$$
\align
&|\Cal G\cap (B_1\times A_2')|> \frac\delta 4|B_1|N_2'\tag 3.9\\
&|\Cal G\cap (A_1'\times B_2)|>\frac \delta 4 |B_2|N_1'\tag 3.10
\endalign
$$
and
$$
N_i'> \frac{3\delta} 4 N_i .\tag 3.11
$$
Moreover, the property
$$
(\Cal G\cap (A_1'\times A_2')^c|\leq \frac\delta 4 |(A_1\times A_2)\backslash (A_1'\times A_2')|\tag 3.12
$$
will hold.
\endproclaim

Thus (3.12) implies (3.11), because
$$
N_1'N_2'\geq |\Cal G\cap(A_1'\times A_2')|> \delta N_1N_2
-\frac\delta 4 N_1N_2=\frac{3\delta}4 N_1N_2.
$$
The construction is straightforward. Assume $A_1'\times A_2'$
fails (3.9). Thus $|\Cal G\cap (B_1\times A_2')|\leq\frac \delta 4
|B_1| \ |A_2'|$ for some $B_1\subset A_1'$. Define $A_1''
=A_1'\backslash B_1$, then
$$
\align
|\Cal G\cap (A_1''\times A_2')^c|&= |\Cal G\cap (A_1'\times A_2')^c|+|\Cal G\cap (B_1\times A_2')|\\
&\leq \frac \delta 4 |(A_1\times A_2)\backslash (A_1'\times A_2')|+\frac \delta 4 |B_1| \ |A_2'|\\
&= \frac\delta 4|(A_1\times A_2)\backslash (A_1'' \times A_2')|
\endalign
$$
and (3.12) remains valid.

Continuing removing the bad set $B_i$, (3.12) ensures that the
remaining set is still big enough, and the process gives the
desired result.
\medskip
\proclaim {Step 2} We decompose $\Cal P_0\subset\Cal P$ in
disjoint sets $\Cal P_0 =\Cal P_1\cup\Cal P_2$.
\endproclaim

 The choice of this decomposition will only matter for condition (3.6).

We proceed as follows

Enumerate $\Cal P_0=\{p_1<p_2<\cdots <p_t\}$ which we identify with $\{1, \ldots, t\}$,

For $t'\leq t$, consider the decreasing functions $(i=1, 2)$
$$
n_i(t') =\max_{(x_1, \ldots, x_{t'})\in \Bbb Z^{t'}}|A_i(x_1,
\ldots, x_{t'})|,
$$
where $A_i(x_1, \ldots, x_{t'})=\{(x_{t'+1}, \ldots, x_t)\mid
(x_1, \ldots, x_t)\in A_i\}.$

We take $t'$ such that
$$
\cases
n_1(t')+n_2(t')\geq (N_1N_2)^{1/4}\\
n_1(t'+1)+n_2(t'+1)\leq (N_1N_2)^{1/4}.\endcases\tag 3.13
$$
We assume $n_1(t')\geq n_2(t')$, thus
$$
n_1(t')\geq \frac 12 N^{1/4}.\tag 3.14
$$
Decompose then $\Cal P_0=\Cal P_1\cup \Cal P_2$ where $\Cal P_1=\{p_1, \ldots, p_{t'}\}$.
\medskip
\proclaim {Step 3} Let $\Cal R_i=\prod_{p\in\Cal P_i}\Bbb Z_{\geq
0} $ with $\Cal R_0=\Cal R_1\times\Cal R_2$ corresponding to the
decomposition in Step 2, and let $\pi_1: \Cal R_0\rightarrow \Cal
R_1$ be the projection to the first $t'$ coordinates. Denote
$$\bar
x=(x_1, \cdots, x_{t'}).$$
 We construct a set ${\bar{\kern-3pt\bar
A}}_2 \subset A_2'$ such that for all $\bar x \in\pi_1(
\,{\bar{\kern-3pt\bar A}}_2),$ we have $|\ {\bar {\kern-3pt\bar
A}}_2(\bar x)|\sim m_2>c\delta^5 K^{-2}N^{1/4}$ and $M_2 \equiv
|\pi_1(\,{\bar{\kern-3pt\bar A}}_2)| \sim \frac {|\
{\bar{\kern-3pt\bar A}}_2|}{m_2}<C\delta^{-5}
K^2\frac{N_2}{N^{1/4}}$, and ${\bar{\kern-2pt\bar N}}_2 \equiv
|\,\bar{\kern-3pt\bar A}_2|> c\frac{\delta^3}{\log\frac
K\delta}N_2.$
\endproclaim

Choose $\bar x\in \pi_1(A_1')$ such that
$$
|A_1'(\bar x)|= n_1(t').\tag 3.15
$$
It follows from (3.9) that
$$
|\Cal G\cap [\big(\{\bar x\}\times A_1'(\bar x)\big)\times
A_2']|>\frac \delta 4 n_1(t')N_2'
$$
and hence there is a subset $A_2''\subset A_2'$ such that by Fact
1 below,
$$
N_2'' =|A_2''|> \frac\delta 8 N_2'\tag 3.16
$$
and for $z\in A_2''$
$$
|\Cal G\cap [\big(\{\bar x\}\times A_1'(\bar x)\big)\times
\{z\}])>\frac\delta 8 n_1(t').\tag 3.17
$$
\medskip
\proclaim {Fact 1} Let $|E|\leq e$ and $|F| \leq f$. If $|\Cal G
\cap (E\times F)| > \alpha ef$, then there exists $F' \subset F$
with $|F'|> \frac {\alpha}{2}f$, such that for any $z \in F'$,
$|\Cal G \cap (E \times \{z\})|> \frac {\alpha}{2} e$.
\endproclaim

 From (2.5) and (3.17), we get clearly
$$
\align
\frac {K^2}\delta N_2\geq K\sqrt{N_1N_2}&= |A_1\underset{\Cal G}\to+ A_2|\\
&\geq |\big(\{\bar x\}\times A_1'(\bar x)\big)\underset {\Cal G}
\to+ A_2''|> \frac \delta 8 |\pi_1(A_2'')|\cdot n_1(t').\tag 3.18
\endalign
$$
Let $\bar A_2 \subset A_2''$ such that the fibers over $\bar x$
have size at least $\frac {\delta^5 n_1(t')}{10^4 K^2}$, thus
$$
\bar A_2 =\bigcup_{|A_2''(\bar x)|> 10^{-4} \delta^5K^{-2}
n_1(t')} \big(\{\bar x\}\times A_2''(\bar x)\big).
$$
It follows from (3.18) that
$$
|A_2''\backslash \bar A_2|\leq |\pi_1(A_2'')| 10^{-4}\delta^5 K^{-2} n_1(t')< \delta^3 10^{-3}N_2 <\frac
\delta{10}N_2''\tag 3.19
$$
The last inequality is by (3.11) and (3.16).

Since by (3.10)
$$
|\Cal G\cap (A_1'\times A_2'')|>\frac \delta 4 N_1' N_2'',
$$
it follows from (3.19) that
$$
|\Cal G\cap(A_1'\times\bar A_2)|>\frac \delta4  N_1' N_2'' -\frac \delta{10}  N_1' N_2''
>\frac\delta{10}  N_1' N_2''.\tag 3.20
$$
Since $|A_2''(\bar x)|\leq n_2(t') \leq n_1(t')$, we may specify
$m_2$ and $\,{\bar{\kern-3pt \bar A}}_2$ as follows:
$$
10^{-4} \delta^5 K^{-2} n_1(t')<m_2<n_1(t'),\tag 3.21
$$
and
$$
A_2' \supset A_2'' \supset \bar A_2 \supset {\bar{\kern-3pt \bar
A}}_2 =\bigcup_{|A_2''(\bar x)|\sim m_2} \big(\{\bar x\}\times
A_2''(\bar x)\big) \tag 3.22
$$
such that
$$
|\Cal G\cap (A_1'\times \, {\bar{\kern-3pt\bar A}}_2)|>c\frac\delta{\log\frac K\delta} N_1'N_2''.\tag 3.23
$$
Thus $\, {\bar{\kern-2pt\bar N}}_2 =|\ {\bar{\kern-3pt\bar A}}_2|$ satisfies
$$
{\bar{\kern-2pt\bar N}}_2 >c\frac\delta {\log\frac K\delta} N_2''
> c\frac{\delta^3}{\log\frac K\delta}N_2.\tag 3.24
$$
The set \ ${\bar{\kern-3pt\bar A}}_2$ has a `regular' structure
with respect to the decomposition $\Cal R_0=\Cal R_1\times\Cal
R_2$ in the sense that for all $\bar x\in\pi_1(
\,{\bar{\kern-3pt\bar A}}_2),$ with $|\ {\bar {\kern-3pt\bar
A}}_2(\bar x)|\sim m_2$. In particular, denoting $M_2=|\pi_1(\,
{\bar {\kern-3pt\bar A}}_2)|$, we have
$$
{\bar{\kern-2pt\bar N}}_2 \sim M_2\cdot m_2.\tag 3.24'
$$
By (3.21) and (3.14)
$$
m_2>c\delta^5 K^{-2}N^{1/4}, $$ and
 $$
M_2 <C\delta^{-5}
K^2\frac{N_2}{N^{1/4}}.\tag 3.25
$$
\medskip
\proclaim {Step 4} Regularization of $A_1'$. We construct a set
$\,\bar{\kern-3pt\bar A}_1 \subset A_1'$ such that for any $\bar x
\in \pi _1 (\,\bar{\kern-3pt\bar A}_1)$, we have
$|\bar{\kern-3pt\bar A}_1(\bar x)|\sim m_1 >c\delta^{10}
K^{-5}N^{1/4},M_1 \equiv |\pi _1 (\,\bar{\kern-3pt\bar A}_1)|\sim
\frac {|\,\bar{\kern-3pt\bar A}_1|}{m_1}< C\delta^{-10}K^5 \frac
{N_1}{N^{1/4}},$ \noindent ${\,\bar{\kern-2pt\bar N}}_1
=|\,\bar{\kern-3pt\bar A}_1|> c\frac{\delta^2}{(\log \frac
K\delta)^2} N_1,$ and $ |\,\Cal G \cap ({\,\bar{\kern-3pt\bar
A}}_1\times{\,\bar{\kern-3pt\bar A}}_2)|> c\frac\delta{(\log \frac
K\delta)^2} {\,\bar {\kern-2pt\bar N}}_1 {\, \bar{\kern-2pt\bar
N}}_2.$
\endproclaim

\proclaim{Claim} Let $\tilde A_1\subset A_1'$. If
$$
|\Cal G\cap (\tilde A_1\times {\,\bar{\kern-3pt\bar A}}_2)|\sim
|\Cal G\cap (A_1' \times{\, \bar{\kern-3pt\bar A}}_2)|, \tag 3.26
$$
then
$$
m\equiv\max _{\bar x \in \pi _1(\tilde A_1)} |\tilde A_1(\bar
x)|\,
> \, c\frac{\delta^4}{(\log\frac K\delta)^2} K^{-2} m_2. \tag 3.27
$$
\endproclaim
\noindent
 {\bf Proof of Claim.} From (3.26), (3.23) and the
regular structure of ${\, \bar{\kern-3pt\bar A}}_2$, there is
$\bar x\in\pi_1( {\, \bar{\kern-3pt\bar A}}_2)$ such that
$$
|\Cal G\cap \big(\tilde A_1\times \big(\{\bar x\}\times{\,
\bar{\kern-3pt\bar A}}_2(\bar x)\big)\big)|
>c\frac\delta{\log\frac K\delta} N_1'm_2.
$$
Hence by Fact 1, there is a subset $A_1'' \subset \tilde A_1$
satisfying
$$
|A_1''|> c\frac\delta{\log\frac K\delta}N_1'\tag 3.28
$$
and for any $z\in A_1''$
$$
|\Cal G \cap (\{z\}\times(\{\bar x\}\times{\, \bar{\kern-3pt\bar
A}}_2(\bar x)\big)\big)|> c\frac \delta{\log\frac K\delta} m_2.
$$
As in Step 3, (3.18), write
$$
\align
\frac {K^2}\delta N_1 \geq K\sqrt{N_1N_2}\geq |A_1\underset{\Cal G}\to+ A_2|&\geq |A_1'' \underset{\Cal G}\to+ \big
(\{\bar x\}\times {\, \bar{\kern-3pt\bar A}}_2 (\bar x)\big)|\\
&>c|\pi_1(A_1'')|\frac\delta{\log\frac K\delta}m_2\\
&> c\frac{|A_1''|}m \ \frac\delta{\log\frac K\delta}m_2\\
&>c\frac {\delta^3}{(\log \frac K\delta)^2} \ \frac{m_2}{m} N_1.
\endalign
$$
The last two inequalities are by the definition of $m$ in (3.27)
and (3.28), (3.11). Hence
$$
m>c\frac{\delta^4}{(\log\frac K\delta)^2} K^{-2} m_2.\tag 3.29
$$
\medskip
Clearly, the bound in (3.27) is smaller than $\delta^5K^{-3}m_2$.
Therefore, in (3.23) we may replace $A_1'$ by $\bar A_1$ defined
as follows.
$$
A_1' \supset\bar A_1=\bigcup_{|A_1'(\bar x)|> \delta^5 K^{-3}m_2}
\big(\{\bar x\}\times A_1'(\bar x)\big).
$$
Thus
$$
|\Cal G\cap (\bar A_1\times{\,\bar{\kern-3pt\bar A}}_2)|>c\frac\delta{\log\frac K\delta} N_1' N_2''.
$$
Recalling (3.21), for $\bar x\in\pi_1(\bar A_1)$
$$
\delta^5K^{-3}m_2< |\bar A_1(\bar x)|\leq n_1(t')<
C\delta^{-5}K^2m_2.
$$
Keeping (3.23) and (3.26) in mind, we may thus again specify
$$
\delta^5K^{-3}m_2 < m_1< C\delta^{-5}K^2 m_2\tag 3.30
$$
such that the regular set $\bar{\kern-3pt\bar A}_1$ defined as
$$
A_1' \supset \bar A_1 \supset{\,\bar{\kern-3pt\bar A}}_1
=\bigcup_{|\bar A_1(\bar x)|\sim m_1} \big(\{\bar x\}\times\bar
A_1(\bar x)\big) \tag 3.31
$$
will satisfy
$$
|\Cal G\cap({\, \bar{\kern-3pt\bar A}}_1\times{\, \bar{\kern-3pt\bar A}}_2)|> c\frac \delta{(\log \frac K\delta)^2} N_1' N_2''.\tag 3.32
$$
Denoting ${\,\bar{\kern-2pt\bar N}}_1 = |{\ \bar{\kern-3pt\bar
A}}_1|, M_1=|\pi_1({\,\bar{\kern-3pt\bar A}}_1)|$, we have
${\,\bar{\kern-3pt\bar N}}_1 \sim M_1m_1$. On the other hand,
(3.32), (3.11) and the fact that $\,\bar {\kern-3pt\bar A}_i
\subset A_i''$ give
$$
{\,\bar{\kern-2pt\bar N}}_1 > c\frac{\delta^2}{(\log \frac K\delta)^2} N_1\tag 3.33
$$
and

$$
|\Cal G \cap ({\,\bar{\kern-3pt\bar
A}}_1\times{\,\bar{\kern-3pt\bar A}}_2)|> c\frac\delta{(\log \frac
K\delta)^2} {\,\bar {\kern-2pt\bar N}}_1 {\, \bar{\kern-2pt\bar
N}}_2. \tag 3.34
$$
It follows from (3.25) and (3.30) that
$$
m_1> c\delta^{10} K^{-5}N^{1/4}, $$ $$
 M_1< C\delta^{-10}K^5 \frac
{N_1}{N^{1/4}}.\tag 3.35
$$
Thus at this stage we have regularized both $A_1, A_2$ with
respect to the decomposition $\Cal P_0=\Cal P_1\cup \Cal P_2$.
\medskip
\noindent {\it For simplicity, we re-denote ${
\bar{\kern-3pt\bar A}}_1, {\, \bar {\kern-3pt\bar A}}_2$ by $A_1,
A_2$ whose cardinalities ${\,\bar{\kern-2pt\bar N}}_i \sim m_i
M_i$ satisfying (3.24) and (3.33).}
\medskip
\proclaim {Step 5} Regularization of the graph. We construct
$\,\Cal G_{1, 0}\subset\pi_1(A_1)\times \pi_1(A_2) \subset \Cal
R_1 \times \Cal R_1$ with $|\,\Cal G_{1,0}|>\delta_0M_1M_2$, such
that $\forall (\bar x_1, \bar x_2)\in \Cal G_{1,0}$, we have $
|A_1(\bar x_1){{\underset{\Cal G}\to+}_{\bar x_1, \bar x_2}}
A_2(\bar x_2)| \sim L\sqrt{m_1m_2},$ and $|\,\Cal G_{\bar x_1,
\bar x_2}|\sim \delta_1 m_1m_2$, where $\Cal G _{\bar x_1, \bar
x_2} $ is the fiber over $(\bar x_1, \bar x_2)$, and $\delta_0,
\delta_1$ and $L$ satisfy (3.45), (3.38) and (3.43) respectively.
\endproclaim

For $\bar x_1, \bar x_2\in \Cal R_1$, let $\Cal G_{\bar x_1, \bar
x_2}$ be the fiber over $(\bar x_1, \bar x_2)$,
$$
\Cal G_{\bar x_1, \bar x_2} =\{(\bar y_1, \bar y_2)\in A_1(\bar
x_1)\times A_2(\bar x_2)| \big((\bar x_1, \bar y_1), (\bar x_2,
\bar y_2)\big)\in\Cal G\}\subset \Cal R_2\times \Cal R_2.
$$
It follows from (3.34) that we may restrict $\Cal G$ to $\Cal
G_1\times (\Cal R_2\times\Cal R_2)$, where
$$
\Cal G_1 =\{(\bar x_1, \bar x_2)\in \pi_1(A_1)\times \pi_1(A_2)| \
|\Cal G_{\bar x_1, \bar x_2}|>c\frac\delta{(\log \frac K\delta)^2}
m_1m_2 \}.\tag 3.36
$$
Thus
$$
c\,m_1m_2\geq |\Cal G_{\bar x_1, \bar x_2}|> c\frac\delta{(\log
\frac K\delta)^2} m_1m_2, \, \text { for } (\bar x_1, \bar x_2)\in
\Cal G_1
$$
and by (3.34)
$$
\sum_{(\bar x_1, \bar x_2)\in\Cal G_1}|\Cal G_{\bar x_1, \bar
x_2}|> c\frac\delta{(\log\frac K\delta)^2} \ {\bar{\kern-2pt\bar
N}}_1 {\, \bar{\kern-2pt\bar N}}_2.\tag 3.37
$$
We may thus specify $\delta_1$,
$$
1>\delta_1> c\frac\delta{(\log\frac K\delta)^2}\tag 3.38
$$
such that if
$$
\Cal G_1' =\{(\bar x_1, \bar x_2)\in\Cal G_1| \ |\Cal G_{\bar x_1,
\bar x_2}|\sim \delta_1 m_1m_2\}, \tag 3.39
$$
then we have
$$
\sum_{(\bar x_1, \bar x_2)\in\Cal G_1'} |\Cal G_{\bar x_1, \bar
x_2}|> c \frac\delta{(\log\frac K\delta)^3}\,
{\,\bar{\kern-2pt\bar N}}_1{\,\bar{\kern-2pt\bar N}}_2. \tag 3.40
$$

Hence
$$
|\Cal G_1'|> c\frac\delta{\delta_1(\log\frac K\delta)^3}
M_1M_2,\tag 3.41
$$
which is bigger than $\frac\delta{(\log\frac K\delta)^3} M_1M_2$.

\noindent
 By further restriction of $\Cal G_1'$, we will also make
a specification on the size of the sumset of $\Cal G_{\bar x_1,
\bar x_2}$.

For $(\bar x_1, \bar x_2) \in \Cal G_1'$, let $K(\Cal G_{\bar x_1,
\bar x_2})$ be the addition constant of $A_1(\bar x_1)$ and
$A_2(\bar x_2)$ along the graph $\Cal G_{\bar x_1, \bar x_2}$ as
defined in Proposition 3. Let $\Cal H \subset \Cal G_1'$, with $$
|\Cal H|\sim |\Cal G_1'|>\frac\delta{(\log\frac K\delta)^3}
M_1M_2.\tag 3.42
$$
\proclaim{Claim} $$ \min_{(\bar x_1, \bar x_2) \in \Cal H} K(\Cal
G_{\bar x_1, \bar x_2}) < L_0 \equiv  (\log \frac
{K}{\delta})^{\frac 92} \delta^{-\frac 92}K.\tag 3.43 $$
\endproclaim

{\bf Proof.} Assume for all $(\bar x_1, \bar x_2) \in \Cal H$
that $K(\Cal G_{\bar x_1, \bar x_2}) > L_0$. Then
$$ \align
K\sqrt{N_1N_2}\geq|A_1\underset{\Cal G}\to+ A_2| &> \min_{(\bar
x_1, \bar x_2) \in \Cal H}\{|A_1(\bar x_1) \, {\underset{\Cal
G}\to +}_{\bar x_1, \bar x_2} A_2(\bar
x_2)|\}|\pi_1(A_1)\underset{\Cal H}\to+ \pi_1(A_2)|\\ &\geq
L_0\sqrt{m_1m_2} \ \frac{|\Cal H|}{\sqrt {M_1M_2}} > L_0\frac {c\,
\delta}{(\log \frac K\delta)^3} \, ({\,\bar{\kern-2pt\bar N}}_1
{\,\bar{\kern-2pt\bar N}}_2)^{1/2}\\ &> c\,\delta^{-1}
\sqrt{N_1N_2}K,
\endalign
$$
which is a contradiction. (The last inequality is by (3.24),
(3.33) and (3.43).)
\medskip

Hence, we may reduce $\Cal G_1'$ to $\Cal G_1''\subset\Cal G_1',$
with $|\Cal G_1''|\sim |\Cal G_1'|$ such that
$$
|A_1(\bar x_1){{\underset{\Cal G}\to+}_{\bar x_1, \bar x_2}}
A_2(\bar x_2)|< L_0\sqrt{m_1m_2}\text { for } (\bar x_1, \bar
x_2)\in \Cal G_1''.
$$
Therefore there is $\Cal G_{1, 0}\subset \Cal G_1''$ and $L<L_0$
(see (3.43))
$$
|\Cal G_{1, 0}|> \frac{|\Cal G_1''|}{\log\frac K\delta}>\delta_0
M_1M_2,\tag 3.44 $$ where, by (3.41)
$$ \delta_0> c \, \frac\delta{
\delta_1(\log\frac K\delta)^4}\tag 3.45
$$
and
$$
|A_1(\bar x_1){{\underset{\Cal G}\to+}_{\bar x_1, \bar x_2}}
A_2(\bar x_2)| \sim L\sqrt{m_1m_2}
\tag 3.46
$$
for $(\bar x_1, \bar x_2)\in\Cal G_{1, 0}.$

Since
$$\align
K\sqrt{N_1N_2}&\geq |\pi_1 (A_1){\underset{\Cal G}\to+}_{1, 0}
\pi_1(A_2)||A_1(\bar x_1){\underset{\Cal G}\to+}_{\bar x_1, \bar
x_2}A_2(\bar x_2)|\\&\geq |\pi_1 (A_1){\underset{\Cal G}\to+}_{1,
0} \pi_1(A_2)|\cdot L\sqrt{m_1m_2}\\&= K(\Cal G_{1, 0})L\sqrt {{\,
\bar{\kern-2pt\bar N}}_1{\,\bar{\kern-2pt\bar N}}_2},
\endalign
$$
we have
$$
K(\Cal G_{1, 0})\cdot L<\delta^{-\frac {5}{2}} (\log \frac {
K}{\delta})^{\frac {3}{2}}K<\delta^{-3} (\log K)^2K.\tag 3.47
$$
In summary, $\Cal G_{1, 0}\subset\pi_1(A_1)\times \pi_1(A_2)$
satisfies (3.44), (3.45) and for $(\bar x_1, \bar x_2)\in\Cal
G_{1, 0}$, the graph $\Cal G_{\bar x_1, \bar x_2}\subset A_1(\bar
x_1)\times A_2(\bar x_2)$ satisfies
$$
\align
&\{(\bar x_1, \bar x_2)\}\times\Cal G_{\bar x_1, \bar x_2}\subset\Cal G\\
&|\Cal G_{\bar x_1, \bar x_2}|\sim \delta_1 m_1m_2,\tag 3.48
\endalign
$$
where $\delta_1$ is as in (3.38). The addition constants $K(\Cal
G_{1, 0})$ and $L$ satisfy (3.43) and (3.47).
\medskip
Denote
$$
\Cal G\supset\tilde{\Cal G} =\bigcup_{(\bar x_1, \bar x_2)\in\Cal
G_{1, 0}} (\{(\bar x_1, \bar x_2)\}\times\Cal G_{\bar x_1, \bar
x_2})\tag 3.49
$$
which satisfies
$$
|\tilde{\Cal G}|> c\, \frac\delta{(\log\frac K\delta)^4} \ {\,\bar{\kern-2pt\bar N}}_1{\,\bar{\kern-2pt\bar N}}_2\tag 3.50
$$
where
$$
{\,\bar{\kern-2pt\bar N}}_1\cdot{\,\bar{\kern-2pt\bar N}}_2> \frac {\delta^5}{(\log \frac K\delta)^3}N_1N_2.\tag 3.51
$$
\medskip
\proclaim {Step 6} Moment inequalities
\endproclaim

Let $\tilde{\Cal G}$ be the graph obtained in (3.49). We reduce
further $\tilde{\Cal G}$ to a graph $\Cal G'$ to fulfil condition
(2.4).

Consider first the graph $\Cal G_{1, 0} \subset\pi_1(A_1)\times\pi_1(A_2)\subset\Cal R_1\times\Cal R_1$ and denote
$K_0=K(\Cal G_{1, 0})$.
By (3.44) and since $\phi, \psi$ are admissible, there is $\Cal G_{1, 0}'\subset\Cal G_{1, 0}$ satisfying
$$
|\Cal G_{1, 0}'|>\phi(M_1 M_2, \delta_0, K_0)\tag 3.52
$$
and
$$
\bigg\Vert\sum F_{\alpha} (\prod_{p\in\Cal P_1} p^{\alpha_p}\theta)\bigg\Vert_q\leq \psi(M_1M_2, \delta_0, K_0)
\bigg(\sum\Vert F_\alpha\Vert_q^2\bigg)^{1/2}\tag 3.53
$$
whenever $\bar x\in\Cal R_1, \, (F_\alpha)_{\alpha\in\Cal G_{1,
0}'(\bar x)}$ trigonometric polynomials satisfying
$$
(n, p)=1, \text { for any } n\in\text{\,supp\,} \widehat F_\alpha,
\text { and for any } p\in\Cal P_1.\tag *
$$
Next, fix $(\bar x_1, \bar x_2)\in\Cal G_{1, 0}'$ and consider the
graph $\Cal G_{\bar x_1, \bar x_2}\subset A_1(\bar x_1)\times
A_2(\bar x_2)\subset \Cal R_2\times\Cal R_2$ satisfying (3.48) and
(3.46).

Thus there is a subgraph $\Cal G_{\bar x_1, \bar x_2}'\subset\Cal
G_{\bar x_1, \bar x_2}$ s.t.
$$
|\Cal G_{\bar x_1, \bar x_2}'| >\phi(m_1m_2, \delta_1, L)\tag 3.54
$$
and
$$
\bigg\Vert\sum G_\alpha\bigg(\prod _{p\in\Cal P_2}p^{\alpha_p}\theta\bigg)\bigg\Vert_q\leq \psi(m_1m_2, \delta_1, L)
\bigg(\sum\Vert G_\alpha\Vert^2_q\bigg)^{1/2}\tag 3.55
$$
whenever $\bar y\in\Cal R_2$ and $(G_\alpha)_{\alpha\in \Cal
G_{\bar x_1, \bar x_2}'(\bar y)}$ trigonometric polynomials
satisfying
$$
(n, p)=1, \text { for any } n\in\text{ supp\,} \widehat G_\alpha,
\text { and for any } p\in\Cal P_2 .\tag {$**$}
$$
Consider then the subgraph $\Cal G'\subset\tilde {\Cal G} \subset \Cal G$
$$
\Cal G'=\bigcup_{(\bar x_1, \bar x_2)\in \Cal G_{1, 0}'} (\{(\bar
x_1, \bar x_2)\}\times \Cal G_{\bar x_1, \bar x_2}')\tag 3.56
$$
which satisfies by (3.52) and (3.54)
$$
|\Cal G'|>\phi(M_1M_2, \delta_0, K_0)\cdot \phi(m_1m_2, \delta_1, L).\tag 3.57
$$
Next, we check the moment inequality.

Fix thus $x=(\bar x, \bar y)\in \Cal R_1\times\Cal R_2$ and
consider trigonometric polynomials $(G_\alpha)_{\alpha\in\Cal
G'(x)}$, where $\Cal G'(x)=\{\alpha \mid (x, \alpha) \in \Cal G'
\subset A \times A \}$, such that
$$
(n, p) =1, \text { for any } n\in\text{\,supp\,} \widehat
G_\alpha, \text { and for any } p\in\Cal P_0=\Cal P_1\cup \Cal
P_2.
$$
We need to estimate $\Vert\sum_\alpha G_\alpha (\prod_{p\in\Cal
P_0} p^{\alpha_p}\theta) \Vert_q$ using (3.53) and (3.55). First,
by (3.56)
$$
\Cal G'(x)=\bigcup_{\beta\in \Cal G_{1, 0}' (\bar x)}
\big(\{\beta\} \times\Cal G_{\bar x, \beta}'(\bar y)\big)
\subset\Cal R_0.
$$
Denote for $\beta\in\Cal G_{1, 0}' (\bar x)$
$$
F_\beta(\theta) =\sum_{\pi_1(\alpha)=\beta} G_\alpha \bigg(\prod_{p\in\Cal P_2} p^{\alpha_p}\theta\bigg)
$$
which clearly satisfy ($*$).

Hence, applying consecutively (3.55) and (3.53)
$$
\align
\bigg\Vert\sum_\alpha G_\alpha \bigg(\prod_{p\in\Cal P_0} p^{\alpha_p}\theta\bigg)\bigg\Vert_q &=\bigg\Vert
\sum_{\beta\in\Cal G_{1, 0}' (\bar x)} F_\beta\bigg(\prod_{p\in\Cal P_1} p^{\beta_p}\theta\bigg)\bigg\Vert_q\\
&\leq\psi (M_1M_2, \delta_0, K_0)\bigg(\sum \Vert F_\beta\Vert^2_q\bigg)^{1/2}\\
&\leq\psi(M_1M_2, \delta_0, K_0).\,\psi (m_1m_2, \delta_1, L)\bigg(\sum_\alpha\Vert
G_\alpha\Vert^2_q\bigg)^{1/2}.\tag 3.58
\endalign
$$
Returning to the statement in Lemma 3.2 and inequalities (3.57)
and (3.58), we get in both (3.3) and (3.4)
$$
\cases
N'=M_1M_2 \qquad N''=m_1m_2\\
\delta'=\delta_0\qquad \delta'' =\delta_1\\
K'=K_0 \qquad K''=L
\endcases
$$
Condition (3.5) follows from (3.51) (which is clearly much
stronger,) and, restating (3.45) and (3.47)
$$
\align
 \delta_0\delta_1 &>\frac\delta{(\log\frac K\delta)^4},\\ K_0
L&<\delta^{-3} (\log K)^2K.\tag 3.59
\endalign
$$
It remains to consider condition (3.6).

By (3.25) and (3.35)
$$
M_1M_2 < C\delta^{-15} K^7 N^{1/2}\tag 3.60
$$
but we don't have necessarily the desired bound on $m_1 m_2$. To
achieve this, we will redefine $\Cal G_{\bar x_1, \bar x_2}'$ by
performing one more step in the construction
\medskip
{\bf Step 7} Recalling Step 2, decompose $\Cal P_2=\{p_{t'+1},
\cdots, p_t\}$ further as
$$
\Cal P_2= \{p_{t'+1}\} \cup \Cal P_3, \text { where } \Cal P_3
=\{p_{t'+2}, \cdots, p_t\}.
$$
For fixed $(\bar x_1, \bar x_2)\in\Cal G_{1, 0}'$, consider the
graph $\Cal K=\Cal G_{\bar x_1, \bar x_2} \subset A_1(\bar
x_1)\times A_2(\bar x_2)\subset \Cal R_2\times\Cal R_2$ satisfying
by (3.46) and (3.48)
$$
\align
&|A_1(\bar x_1)|\sim m_1, \, \, \,|A_2(\bar x_2)|\sim m_2\\
&|\Cal G_{\bar x_1, \bar x_2}|\sim \delta_1 m_1m_2\\
&K(\Cal G_{\bar x_1, \bar x_2})\sim L.
\endalign
$$
Repeat then Steps 1 - 5 from previous construction to the graph
$\Cal K$ with respect to the decomposition $\Cal
P_2=\{p_{t'+1}\}\cup \Cal P_3$. Thus $\Cal K$ gets replaced by
$$
\tilde{\Cal K} =\bigcup_{(z_1, z_2)\in \Cal K_{1, 0}}\Cal K_{z_1, z_2}\tag 3.61
$$
where
$$
\Cal K_{1, 0} \subset (\Bbb Z_{\geq 0})^2.
$$
Thus $\tilde{\Cal K} \subset{\,\bar{\kern-3pt\bar A}}_1(\bar
x_1)\times {\,\bar{\kern-3pt\bar A}}_2 (\bar x_2) \subset A_1(\bar
x_1)\times A_2(\bar x_2)$,
$$
\align
&\Cal K_{z_1, z_2}\subset{\,\bar{\kern-3pt\bar A}}_1(\bar x_1, z_1)\times {\ \bar{\kern-3pt\bar A}}_2(\bar x_2, z_2)\\
&m_i\geq |{\ \bar{\kern-3pt\bar A}}_i (\bar x_i)|\equiv
{\,\bar{\kern-1pt\bar m}}_i>
\frac{\delta_1^3}{(\log\frac L{\delta_1})^2} m_i\tag 3.62\\
&|{\,\bar{\kern-3pt\bar A}}_i(\bar x_i, z_i)|\sim\ell_i \leq
|A_i(\bar x_i, z_i)|< (N_1N_2)^{1/4}\tag 3.63
\endalign
$$
\big(by (3.13)\big)
$$
\align
|\Cal K_{z_1, z_2}|&\sim \delta_3 \ell_1 \ell_2\tag 3.64\\
|\Cal K_{1, 0}|&> \frac{\delta_1}{\delta_3(\log\frac
L{\delta_1})^4} \ \frac {{\, \bar{\kern-.75pt\bar m}}_1
{\,\bar{\kern-.75pt\bar m}}_2} {\ell_1 \ell_2} \tag 3.65
\endalign
$$
(cf. (3.44), (3.45))
$$
K(\Cal K_{z_1, z_2})<K(\Cal K_{1,0})\cdot K(\Cal K_{z_1, z_2})<
\delta_1^{-3} (\log L)^2L\tag 3.66
$$
\big(cf. (3.47)\big).

(We point out here that $\ell_i, {\,\bar{\kern-1pt\bar m}}_i,
\delta_3> \frac{\delta_1}{(\log \frac L{\delta_1})^2}$ do depend
on the basepoint $(\bar x_1, \bar x_2)\in \Cal R_1\times\Cal
R_1)$. Starting from (3.61), we carry out Step 6. However, since
$\Cal K_{1, 0}\subset(\Bbb Z_{\geq 0})^2$ (only the prime
$p_{t'+1}$ is involved) we may take $\Cal K_{1, 0}'=\Cal K_{1, 0}$
and replace in (3.53) the factor $\psi( \ )$ by $Cq$ (apply
Proposition 1 with $k=1$). For each $(z_1, z_2)\in \Cal K_{1, 0}$,
consider again a subgraph $\Cal K_{z_1, z_2}'\subset\Cal K_{z_1,
z_2}$ satisfying
$$
|\Cal K_{z_1,z_2}'| >\phi\big(\ell_1\ell_2, \delta_3, K(\Cal
K_{z_1, z_2})\big)\tag 3.67
$$
and
$$
\bigg\Vert\sum G_\alpha\bigg(\prod_{p\in\Cal P_3} p^{\alpha_p}\theta\bigg)\bigg\Vert_q=\psi \big(\ell_1\ell_2, \delta_3,
K(\Cal K_{z_1, z_2})\big) \bigg(\sum\Vert G_\alpha\Vert^2_q\bigg)^{1/2}\tag 3.68
$$
whenever $(G_\alpha)_{\alpha\in\Cal K_{z_1, z_2}'} (\bar y)$ are
trigonometric polynomials satisfying
$$
(n, p)=1, \text { for any } n\in \text{\, supp\,} \widehat
G_\alpha, \text { and for any } p\in\Cal P_3.
$$
Redefine then $\Cal G_{\bar x_1, \bar x_2}' \subset\Cal G_{\bar
x_1, \bar x_2}$ as
$$
\Cal G_{\bar x_1, \bar x_2}' =\bigcup_{(z_1, z_2)\in \Cal K_{1,
0}}(\{(z_1, z_2)\} \times\Cal K_{z_1, z_2}')\tag 3.69
$$
and take again
$$
\Cal G' =\bigcup_{(\bar x_1, \bar x_2)\in\Cal G_{1, 0}'} (\{(\bar
x_1, \bar x_2)\}\times\Cal G_{\bar x_1, \bar x_2}').
$$
>From the preceding, since $\ell_1\ell_2 \leq \min \{m_1m_2,
N^{1/2}\}$ in (3.63) and (3.66), the factor in the moment bound
(3.59) becomes now
$$
Cq\psi(M_1M_2,\, \delta_0, \,K_0)\cdot \psi\bigg(\min \{m_1m_2,
N^{1/2}\}, \,\frac{\delta_1}{(\log\frac L{\delta_1})^2},\,
\delta_1^{-3}(\log L)^2L\bigg).
$$
Thus in (3.4), $N'=M_1M_2$, $N'' =\min \{m_1m_2 , N^{1/2}\}$
satisfy (3.6) \big(and may be increased to satisfy also the lower
bound in (3.5)\big).

Also
$$
\delta_0\cdot\frac{\delta_1}{(\log \frac L{\delta_1})^2}\overset{(3.51)}\to > c\frac\delta{(\log\frac K\delta)^6}
$$
which is condition (3.7).

Taking $K' =K_0, K''=\delta_1^{-3}(\log L)^2L$, (3.59) implies
$$
K'\cdot K'' <\delta^{-6} (\log K)^{17}K
$$
and hence (3.8) holds.

>From (3.65), (3.67), (3.69)
$$\align
|\Cal G_{\bar x_1, \bar
x_2}'|&>\bigg\{1+\frac{\delta_1}{\delta_3(\log\frac K\delta)^4}
\frac{{\,\bar{\kern-1pt\bar m}}_1
{\,\bar{\kern-1pt\bar m}}_2}{\ell_1 \ell_2}\bigg\}\cdot \phi\bigg(\ell_1\ell_2, \,\delta_3, \,\delta_1^{-3} (\log L)^2L\bigg)\\
&\overset (3.62)\to > \bigg\{1+\frac{\delta_1^7}{(\log \frac
K\delta)^8} \, \frac {m_1m_2}{\ell_1\ell_2}\bigg\}
\cdot\phi\bigg(\ell_1\ell_2, \,\frac{\delta_1}{(\log\frac
K\delta)^2}, \,\delta_1^{-3} (\log L)^2L\bigg). \tag 3.70
\endalign
$$
Define
$$
N'' =\min\bigg \{N^{1/2}, \,1+\frac{\delta_1^7}{(\log \frac
K\delta)^8}m_1m_2\bigg\}.\tag 3.71
$$
Using property (3.1) of the function $\phi$, we verify that
$$
(3.70)>\bigg(1+\frac{\delta_1^7m_1m_2}{(\log \frac
K\delta)^8}\bigg) \frac 1{N''}\cdot\phi \bigg(N'',\,
\frac{\delta_1}{(\log\frac K\delta)^2}, \,\delta_1^{-3}(\log
L)^2L\bigg).
$$

Hence, again by (3.1)
$$
\align |\Cal G'|&> \bigg(1+\frac{\delta_1^7 m_1m_2}{(\log \frac
K\delta)^8}\bigg) \frac 1{N''}\cdot\phi(M_1 M_2, \,\delta_0,\,
K_0)\cdot \phi \bigg(N'',\, \frac{\delta_1}{(\log\frac
K\delta)^2},
\,\delta_1^{-3}(\log L)^2 L\bigg)\\
&>\phi(N', \,\delta_0,\, K_0)\cdot \phi\bigg(N'',
\,\frac{\delta_1}{(\log \frac K\delta)^2}, \,\delta_1^{-3} (\log
L)^2 L\bigg)\tag 3.72
\endalign
$$
denoting
$$
N'=\bigg(1+\frac{\delta_1^7m_1m_2}{(\log\frac K\delta)^8}\bigg)
\frac{M_1M_2}{N''}.\tag 3.73
$$
Thus
$$
N>N'N''>\frac{\delta_1^7}{(\log \frac
K\delta)^8}\cdot{\,\bar{\kern-2pt\bar N}}_1 {\,\bar{\kern-2pt\bar
N}}_2 > \delta^{12}\bigg(\log\frac K\delta\bigg)^{-25}
N\Rightarrow (3.5)
$$
and by (3.71)
$$
N''\leq N^{1/2}, N'< M_1M_2+\frac{{\,\bar{\kern-2pt\bar
N}}_1{\,\bar{\kern-2pt\bar N}}_2}{N^{1/2}} \overset{(3.60)}\to <
C\delta^{-15} K^7N^{1/2}\Rightarrow (3.6).
$$
This proves Lemma 3.2.

\bigskip
\noindent {\bf Section 4. Proof of Proposition 3, Part II}

Recalling (2.13) and (2.14), we start from the pair of admissible
functions
$$
\align
\phi(N, \delta, K)&= \bigg(\frac\delta K\bigg)^CN\tag 4.1\\
\psi (N, \delta, k)&= \min (q^{(\frac K\delta)^C},\, N^{1/2})\tag
4.2
\endalign
$$
($C$ = some constant).
The $N^{\frac 12}$-bound in (4.2) is obtained from the obvious estimate
$$
\bigg\Vert\sum_\alpha F_\alpha \bigg(\prod_{p\in\Cal P_0}p^{\alpha_p}\theta\bigg)\bigg\Vert_q \leq
\sum_\alpha \Vert F_\alpha\Vert_q \leq N^{1/2}\bigg(\sum\Vert F_\alpha\Vert^2_q\big)^{1/2}
\tag 4.3
$$
(since $\alpha$ ranges in a set of size at most $N$).

Starting from (4.1), (4.2), we produce here a new pair of
admissible functions by applying Lemma 3.2. The next statement
does not yet imply Proposition 2 but displays already a much
better behavior of $\psi$ in $K$.

\proclaim
{Lemma 4.3} Take
$$
\tilde\phi(N, \delta, K)= \bigg(\frac \delta K\bigg)^{C\log\log\frac K\delta}\cdot N\tag 4.4
$$
and
$$
\tilde\psi(N, \delta, K)=q^{(\log\frac K\delta)^{C/\gamma}}\cdot N^\gamma
\tag 4.5
$$
with $C$ an appropriate constant and $0<\gamma<1$ arbitrary.

Then $\tilde\phi, \tilde\psi$ are admissible.
\endproclaim

\noindent
{\bf Proof.}

We will make an iterated application of Lemma 3.2.

Fix $N,\delta, K$ and choose an integer $t$ of the form $2^\ell$ (to be specified).
Starting from $\phi_0=\phi, \psi_0=\psi$, define recursively for $\ell'=0, 1, \ldots,
\ell-1$
$$
\align
\phi_{\ell'+1}(N, \delta, K)&= \min\phi_{\ell'} (N', \delta', K')\cdot\phi_{\ell'} (N'', \delta'', K'')\tag 4.6\\
\psi_{\ell'+1} (N, \delta, K)&= \max Cq\psi_{\ell'}(N',\delta',K')\cdot \psi_{\ell'} (N'',\delta'', K'').\tag 4.7
\endalign
$$
where in (4.6), (4.7) the parameters $N', N'', \delta', \delta'', K', K''$ satisfy (3.4)-(3.7).

We evaluate $\tilde\phi=\phi_\ell, \tilde\psi=\psi_\ell$.

Iterating (4.6), we obtain clearly
$$
\tilde\phi(N, \delta, K)=\prod_{\nu\in\{0, 1\}^\ell}\phi(N_\nu, \delta_\nu, K_\nu)\tag 4.8
$$
where
$$
(N_\nu)_{\nu\in\operatornamewithlimits\cup\limits_{\ell'\leq\ell}\{0,
1\}^{\ell'}}, \quad (\delta_\nu)_{\nu\in
\operatornamewithlimits\cup\limits _{\ell'\leq\ell}\{0,
1\}^{\ell'}}, \quad
(K_\nu)_{\nu\in\operatornamewithlimits\cup\limits_{\ell'\leq
\ell}\{0, 1\}^{\ell'}}
$$
satisfy by (3.4)-(3.7) the following constraints
$$
\align
N_\phi&=N, \delta_\phi=\delta, K_\phi=K\tag 4.9\\
N_\nu&\geq N_{\nu, 0}\cdot N_{\nu, 1}\geq N_\nu\bigg(\frac{\delta_\nu}{\log K_\nu}\bigg)^{40} \tag 4.10\\
N_{\nu, 0}+N_{\nu, 1}&< \bigg(\frac{K_\nu}{\delta_\nu}\bigg)^{20} N_\nu^{1/2}\tag 4.11\\
\delta_{\nu, 0}\cdot\delta_{\nu, 1}&> \bigg(\log\frac{K_\nu}{\delta_\nu}\bigg)^{-6}\delta_\nu.\tag 4.12\\
K_{\nu, 0}\cdot K_{\nu, 1}&\leq \delta_\nu^{-6} (\log K_\nu)^{20} K_\nu.\tag 4.13
\endalign
$$
>From (4.12), (4.13)
$$
\log\frac{K_{\nu, 0}}{\delta_{\nu, 0}} +\log \frac{K_{\nu, 1}}{\delta_{\nu, 1}} < 8\log\frac{K_\nu}{\delta_\nu}
$$
and iteration implies
$$
\max_{\nu \in\{0, 1\}^{\ell'}}\log \frac {K_\nu}{\delta_{\nu}}\leq
\sum_{\nu\in\{0, 1\}^{\ell'}} \log\frac{K_\nu}{\delta_{\nu}}<
8^{\ell'}\log\frac K\delta.\tag 4.14
$$
Iteration of (4.12) gives
$$
\spreadlines{10pt}
\align
\prod_{\nu\in\{0, 1\}^{\ell'}}\delta_\nu &>\prod_{\nu\in\{0, 1\}^{\ell'-1}} \bigg
(\log\frac{K_\nu}{\delta_\nu}\bigg)^{-6} \prod_{\nu\in\{0, 1\}^{\ell'-1}} \delta_\nu\\
&\overset{(4.14)}\to > 8^{-3\ell'2^{\ell'}} \bigg(\log \frac {K}{\delta}\bigg)^{-3\cdot 2^{\ell'}}\prod_{\nu\in\{0, 1\}^{\ell'-1}}\delta_\nu\\
&> 8^{-3(\ell'2^{\ell'}+(\ell'-1)
2^{\ell'-1}+\cdots)}\bigg(\log\frac
{K}{\delta}\bigg)^{-3(2^{\ell'}+2^{\ell'-1}+\cdots)} \delta \\&
>8^{-6\ell' 2^{\ell'}}\bigg(\log\frac{K}{\delta}\bigg)^{-6\cdot
2^{\ell'}}
 \delta.\tag 4.15
\endalign
$$
Next. iterate (4.13).
Thus
$$
\align
\prod_{\nu\in\{0, 1\}^{\ell'}}K_\nu &\leq \prod_{\nu\in\{0, 1\}^{\ell'-1}} \delta_\nu^{-6}(\log K_\nu)^{20}\cdot
\prod_{\nu\in\{0, 1\}^{\ell'-1}}K_\nu\\
&\overset{(4.14), (4.15)}\to < \ \Bigg(8^{-3\ell'2^{\ell'}}
\bigg(\log\frac{K}{\delta}\bigg)^{-3\cdot 2^{\ell'}}
\delta\Bigg)^{-6}\bigg(8^{\ell'}\log\frac
K\delta\bigg)^{10\cdot2^{\ell'}}
\bigg(\prod_{\nu\in \{0, 1\}^{\ell'-1}}K_\nu\bigg)\\
&< 8^{56\ell'2^{\ell'}}\bigg(\log \frac K\delta\bigg)^{56\cdot
2^{\ell'}}\delta^{-6\ell'}K.\tag 4.16
\endalign
$$

By (4.10)
$$
\spreadlines{10pt}
\align
\prod_{\nu\in\{0, 1\}^{\ell'}}N_\nu&> \prod_{\nu\in\{0, 1\}^{\ell'-1}} \bigg(\frac{\delta_\nu}{\log
K_\nu}\bigg)^{40}\cdot \prod_{\nu\in\{0, 1\}^{\ell'-1}}N_\nu\\
&>8^{-140\ell'2^{\ell'}}\delta^{40} \bigg(8^{\ell'}\log\frac
K\delta\bigg)^{-140\cdot 2^{\ell'}} \prod_{\nu\in \{0,
1\}^{\ell'-1}}N_\nu\\
&> 8^{-280\ell' 2^{\ell'}}\bigg(\log\frac K\delta\bigg)^{-280 \,
2^{\ell'}}\delta^{40\ell'}N.\tag 4.17
\endalign
$$

Substitution of (4.1) in (4.8) gives by (4.15), (4.16), (4.17)
$$
\align
\tilde\phi(N, \delta, K)&\geq \prod_{\nu\in\{0, 1\}^\ell} \bigg(\frac{\delta_\nu}{K_\nu}\bigg)^C N_\nu\\
&> e^{-C\ell2^\ell}\bigg(\log\frac K\delta\bigg)^{-C2^\ell}\delta^{C\ell}K^{-C}N\\
&>t^{-Ct}\bigg(\log\frac K\delta\bigg)^{-Ct}(\delta^{C\log t}) K^{-C}N.\tag 4.18
\endalign
$$
Similarly, we will iterate (4.7) with (possibly different)
parameters $(N_\nu), (\delta_\nu), (K_\nu)$ still satisfying
(4.9)-(4.17).

By (4.2)
$$
\tilde\psi(N, \delta, K)=(Cq)^t \, \prod_{\nu\in \{0, 1\}^\ell} \min \big(q^{(\frac{K_\nu}{\delta_\nu})^C},
N_{\nu}^{1/2}\big).\tag 4.19
$$
>From (4.12)(which implies that $\delta_{\nu, 0},\,\delta_{\nu, 1}>
\bigg(\log\frac{K_\nu}{\delta_\nu}\bigg)^{-6}\delta_\nu$) and
(4.14) that
$$
\delta_\nu> 8^{-6\ell^2}\bigg(\log\frac K\delta\bigg)^{-6\ell}\delta\tag 4.20
$$
and from (4.13) (which implies that $K_{\nu, 0}, K_{\nu, 1}\leq
\delta_\nu^{-6} (\log K_\nu)^{20} K_\nu$), (4.14), and (4.20) that
$$
K_\nu < 8^{37\ell^3} \bigg(\log\frac K\delta\bigg)^{37\ell^2}
\delta^{-6\ell}K.\tag 4.21
$$
Hence from (4.11), (4.20), (4.21)
$$
 N_{\nu, 0}+N_{\nu, 1}< 8^{800\ell^3}\bigg(\log\frac
K\delta\bigg)^{800\ell^2}\delta^{-800\ell}K^{20} N_\nu^{1/2}
$$
$$
\align
 N_{\nu}<&\,\bigg(8^{800\ell^3}\big(\log\frac
K\delta\big)^{800\ell^2}\delta^{-800\ell}K^{20}\bigg)^{(1+\frac
12+\frac 14 +\cdot)}N^{\frac{1}{2^{\ell}}}\\<&\,
10^{10^3\ell^3}\bigg(\log\frac K\delta\bigg)^{10^4\ell^2}
\delta^{-10^4\ell}K^{40}N^{\frac 1t} \text { for } \nu\in\{0,
1\}^\ell.\tag 4.22\\
\endalign
$$
To bound (4.19), let $A$ be a number to be specified and partition
$$
\{0, 1\}^\ell =I\cup J\text { with } I=\bigg\{\nu\in \{0, 1\}^\ell\bigg| \frac{K_\nu}{\delta_\nu} \leq A\bigg\}.
$$
It follows then from (4.19) and (4.22) that
$$
(4.19)< q^{A^Ct} \bigg[C^{\ell^3}\bigg(\log\frac
K\delta\bigg)^{10^4\ell^2}\delta^{-10^4\ell}K^{40}N^{\frac
1t}\bigg]^{|J|}.\tag 4.23
$$
Now, we will estimate $|J|$.

>From (4.15), (4.16)
$$
A^{|J|} \leq\prod_{\nu\in\{0, 1\}^\ell} \frac{K_\nu}{\delta_\nu} <
8^{62\ell t}\cdot \bigg(\log\frac K\delta\bigg) ^{62t}
\delta^{-7\ell}K.\tag 4.24
$$
Take
$$
2^\ell =t\sim \log\frac K\delta\tag 4.25
$$
and fixing $0<\gamma<1$, take
$$
\log A\sim \gamma^{-1}\log t.\tag 4.26
$$
It follows then from (4.24) that
$$
|J|<\frac{10^3t\cdot \log t}{\log A} <\gamma t.\tag 4.27
$$
Hence clearly from (4.23) that
$$
\tilde\psi(N, \delta, K)< q^{A^Ct} e^{t^2\log t}N^\gamma < q
^{(\log \frac K\delta)^{C/\gamma}} \cdot N^\gamma\tag 4.28
$$
which is (4.5).

Substitution of (4.25) in (4.18) gives (4.4).

This proves Lemma 4.3.
\bigskip

\noindent {\bf Section 5. Proof of Proposition 3, Part III}

We will first use Lemma 3.2 and Lemma 4.3 to show

\proclaim {Lemma 5.1} Assume again the moment $q$ fixed. Given
$0<\tau,\gamma<\frac 12$, for $i=1,2,3,$ there are positive
constants $A_i=A_i(\tau, \gamma), B_i=B_i(\tau,\gamma)$ such that
taking $N$ sufficiently large
$$
\qquad\cases \phi(N, \delta, K)= K^{-A_1} \, \delta^{A_2\log\log
N} \, e^{A_3(\log\log N)^2} \ N^{1-\tau}\qquad
\qquad  &\text{\rm (5.2)}\\
\psi (N, \delta, K) = K^{B_1} \, \delta^{-B_2\log\log N} \ e^{-B_3(\log\log N)^2} N^\gamma\qquad
\qquad\qquad&\text{\rm (5.3)}
\endcases
$$
is a pair of admissible functions.
\endproclaim

\noindent
{\bf Proof.}

We will proceed in 2 steps.

First, some

\noindent {\bf Notation.} We use `$\ell\ell$' to denote
`$\log\log$'. \medskip
 It follows from Lemma 4.3 (by taking $\frac {\gamma}{4}$ and assuming $\frac {K}{\delta} <N$) that
$$
\qquad\cases
\phi(N, \delta, K)= \big(\frac \delta K\big)^{C_0 \ell\ell N}N\qquad&(5.4)\\
\psi(N, \delta, K)=\min \big\{\exp \big[\log q\cdot\big(\log\frac
K\delta\big) ^{\frac {4^{C_0}}\gamma}\big] \cdot N^{\frac \gamma
4}, \, N^{1/2}\big \}\qquad\qquad&(5.5)
\endcases
$$
are admissible.

First, fix a large integer $\bar N$ (depending on $\tau, \gamma)$ and define
$$
A_3\sim A_2 \sim A_1= C_0 \, \ell\ell{\bar N}\tag 5.6
$$
($A_2, A_3$ will be specified later).

Thus the expression in (5.2) is at most
$$
\bigg(\frac\delta K\bigg)^{C_0\ell \ell{\bar N}} \,
e^{C(\ell\ell{\bar N})^3} N^{1-\tau}< \bigg(\frac\delta
K\bigg)^{C_0\ell\ell N}N=(5.4),
$$
if $\bar N$ is large enough and $N\leq \bar N, \log N\sim\log \bar N$.

Taking $N\leq \bar N,$ then $(5.5) < N^{\frac\gamma 2}<(5.3)$,
provided
$$
\log q\cdot \bigg(\log\frac K\delta\bigg)^{\frac{4C_0}\gamma} < \frac \gamma 4 \log N\tag 5.7
$$
and
$$
e^{B_3(\ell\ell N)^2} < N^{\frac \gamma 2}.\tag 5.8
$$
If (5.7) does not hold, then for some $c>0$,
$$
\frac K\delta> e^{(\log N)^{c\gamma}}.\tag 5.9
$$
Thus if we take
$$
B_3\sim B_2\sim B_1=(\log \bar N)^{1-c\gamma}\tag 5.10
$$
(also $B_2, B_3$ to be specified later). (5.9) implies that
$(5.3)>N$ if $N\leq\bar N$ and (5.8) holds. From (5.10) and this
choice of $B_3$, clearly

 (5.8) holds for $N\leq \bar N$, $\log
N\sim \log \bar N$. Thus with preceding choice of $A_1, A_2, A_3,
B_1, B_2, B_3$, (5.2), (5.3) is admissible in the range $N\leq\bar
N$, $\log N\sim \log\bar N$.

Next we use Lemma 3.2 to establish that (5.2), (5.3) are also
admissible in the range $N>\bar N$. We will proceed by induction
on the size of $N$. Obviously $(5.3)>N^{1/2}$ if $(\frac
K\delta)^{B_1}> N^{1/2}$. Hence we may assume
$$
\frac K\delta< N^{10^{-6}}.\tag 5.11
$$
We want to reduce $N$ by Lemma 3.2. Thus we estimate
$$
\phi(N', \delta', K')\cdot \phi(N'', \delta'', K'')\tag 5.12
$$
from below and
$$
\psi(N', \delta', K') \cdot \psi(N'', \delta'', K'')\tag 5.13
$$
from above, where $N'N'',\delta', \delta'', K', K''$ satisfy
(3.5)-(3.8). Hence for $N>>0$,
$$
N\geq  N'N''> N\bigg(\frac\delta{\log N}\bigg)^{40}\underset{(5.11)}\to> N^{\frac{99}{100}}\tag 5.14
$$
$$
\align
N'+N'' &< \bigg(\frac K\delta\bigg) ^{20}N^{1/2} {\underset(5.11)\to<} N^{\frac{11}{20}}\tag 5.15\\
\delta'\delta''&> \frac {10^{36}\delta}{(\log N)^6}\tag 5.16\\
K'K''&< \delta^{-6} (\log N)^{20}K.\tag 5.17
\endalign
$$
Condition (5.15) reduces indeed $N$ to scale $N^{\frac {11}{20}}$
for which we assume (5.2), (5.3) admissible (notice that since
$N>\bar N$, $\log N', \log N''\gtrsim \log\bar N$). In what
follows, the role of the additional technical factors in (5.2),
(5.3) will become apparent.

Substituting (5.2) in (5.12), we get
$$
(K'K'')^{-A_1} (\delta')^{A_2\ell\ell N'}(\delta'')^{A_2\ell\ell N''} \, e^{A_3[\ell\ell N')^2+(\ell\ell N'')^2]}
(N'N'')^{1-\tau}.\tag 5.18
$$
>From (5.14), (5.15)
$$
\align
N^{\frac{11}{25}} &< N', N''<N^{\frac{11}{20}}\\
\frac{99}{100} \ell\ell N&<\ell\ell N-\log \frac
{25}{11}<\ell\ell
N', \ell\ell N''< \ell\ell N-\log\frac{20}{11}.\tag 5.19
\endalign
$$
>From (5.14)-(5.17), (5.19)
$$
\align (5.18)&> \,\delta^{6A_1}(\log
N)^{-20A_1}K^{-A_1}\bigg[\frac\delta{(\log
N)^6}\bigg]^{A_2(\ell\ell N-\log
\frac{20}{11})} \, e^{\frac{19}{10}A_3(\ell\ell N)^2} N^{1-\tau}\bigg(\frac\delta{\log N}\bigg)^{40(1-\tau)}\\
&> \,K^{-A_1}\cdot\delta^{A_2\ell\ell N}\ e^{A_3(\ell\ell N)^2}
N^{1-\tau}\cdot u \cdot v,
\endalign
$$
where
$$
\align
 u=&\,(\log N)^{-20A_1-6A_2\ell\ell N-40} \ e^{\frac
9{10}A_3(\ell\ell N)^2}\tag 5.20\\
 v=&\,\delta^{6A_1-(\log\frac{20}{11})A_2+40}.\tag 5.21\\
 \endalign
$$
and each of the factors (5.20) and (5.21) will be at least 1 for
suitable choices $A_1<A_2<A_3$ ($A_1\sim A_2\sim A_3$). Hence
(5.18) still admits (5.2) as lower bound.

Similarly, substituting (5.3) in (5.13), we get
$$
\align
&(K'K'')^{B_1}(\delta')^{-B_2\ell\ell N'}(\delta'')^{-B_2\ell\ell N''}\, e^{-B_3[(\ell\ell N')^2+(\ell\ell N'')^2]}
(N'N'')^\gamma\\
<&\,\delta^{-6B_1}(\log N)^{20 B_1} K^{B_1}\bigg[\frac{(\log
N)^6}\delta\bigg]^{B_2(\ell\ell N-\log \frac{20}{11})}
\, e^{-\frac{19}{10} B_3(\ell\ell N)^2} \, N^\gamma \\
<&\,K^{B_1}\cdot \delta^{-B_2\ell\ell N}\, e^{-B_3(\ell\ell
N)^2}N^\gamma\cdot u' \cdot v', \endalign
$$ where
$$
\align u'=&\,(\log N)^{20B_1+6B_2\ell\ell N} \cdot e^{-\frac
9{10}B_3(\ell\ell N)^2}\tag 5.22\\
v'=&\, \delta^{-6B_1+B_2\log\frac{20}{11}}.\tag 5.23\\
\endalign
$$
Again a choice $B_3>B_2>B_1$ ($B_3 \sim B_2\sim B_1$) allows us to
get (5.22) and (5.23) at most 1.

Thus (5.13) still satisfies (5.3).

This proves Lemma 5.1.
\medskip
\noindent {\bf Conclusion of the proof of Proposition 3}

Immediate from Lemma 5.1.

First, we may assume $K=K(\Cal G)< N^{\frac 1\Lambda}$ since (2.2) is otherwise obvious.

Apply Lemma 5.1 with $\tau, \gamma$ replaced by $\frac\tau2, \frac\gamma 2$ and let
$$
\Lambda =\frac{2A_1}{\tau} +A_2+B_1+\frac{2B_2}{\gamma}.\tag 5.24
$$
The choice of $\Lambda$ implies that $\Lambda >\frac
{2A_1}{\tau}$, $\Lambda >B_1$, and $\frac {\Lambda}{2B_2}\gamma
>1$.
 In (5.2)
$$
\phi(N, \delta, K)> \delta^{A_2\ell\ell N}N^{1-\frac\tau 2-\frac{A_1}\Lambda}>\delta^{\Lambda \ell\ell N}
N^{1-\tau}.\tag 5.25
$$
In (5.3)
$$
\psi(N, \delta, K)< K^\Lambda \cdot \delta^{-B_2\ell\ell N}N^{\frac \gamma 2}.\tag 5.26
$$
If $\delta^{-B_2\ell\ell N}< N^{\frac \gamma 2},$ then
$(5.26)<K^\Lambda N^\gamma$. Otherwise in (2.1)
$$
\delta^{\Lambda\ell\ell N}N^{1-\tau}< N^{-\frac{\Lambda}{2B_2}\gamma} N^{1-\tau}< 1,
$$
hence the statement becomes trivial.

\bigskip
\noindent {\bf Section 6. Remarks.}

(1) Going back to (1.12)-(1.15) and the proof of Lemma 5.1 and
Proposition 3, an inspection of the argument shows that one may
take $k(b)=C^{b^4}$ in the Theorem (for some constant $C$). We
certainly did not try to proceed efficiently here.

(2) The proof of the Theorem shows in fact the following stronger statement:

For all $b\in\Bbb Z_+, \delta>0$, there is $k=k(b, \delta)\in \Bbb
Z_+$ such that whenever $A\subset\Bbb Z$, $|A|=N$ sufficiently
large, then either $|kA|>N^b$ or $|A_1^{(k)}|> N^b$ for all
$A_1\subset A, |A_1|>N^\delta$.

(3) As in [Ch], our approach uses strongly prime factorization in $\Bbb Z$.
Thus the argument at this point does not apply to subsets $A\subset\Bbb R$.

\enddocument